\documentclass{commat}

\usepackage{longtable}
\newcommand{\A}{\mathbf{A}}

\title{%
    Local derivations of semisimple Leibniz algebras
}

\author{%
    Ivan Kaygorodov, Karimbergen Kudaybergenov and Inomjon Yuldashev
}

\affiliation{
    \address{Ivan Kaygorodov --
    Centro de Matemática e Aplicações, Universidade da Beira Interior, Covilh\~{a}, Portugal \\
    Siberian Federal University, Krasnoyarsk, Russia
}
    \email{%
    kaygorodov.ivan@gmail.com
}
    \address{Karimbergen Kudaybergenov --
    V. I. Romanovskiy Institute of Mathematics, Uzbekistan Academy of Sciences, Tashkent, Uzbekistan \\
    Karakalpak State University, Nukus, Uzbekistan
}
    \email{%
    karim2006@mail.ru
}
    \address{Inomjon Yuldashev --
    Nukus State Pedagogical Institute, Nukus, Uzbekistan
}
    \email{%
    i.yuldashev1990@mail.ru
}
}

\thanks{%
	This work is supported by the Russian Science Foundation under grant 18-71-10007.
}

\abstract{%
    We prove that every local derivation on a complex semisimple finite-dimensional Leibniz algebra is a derivation.
}

\keywords{%
    Lie algebra, Leibniz algebra, simple algebra, semisimple algebra, local derivation.
}

\msc{%
   	17A32, 17A60, 17B10, 17B20.
}

\VOLUME{30}
\NUMBER{2}
\firstpage{1}
\DOI{https://doi.org/10.46298/cm.9274}

\begin{paper}

The study of local derivations
started with Kadison’s article~\cite{Kadison90}.
After this work, appear numerous new results related to the description of numerous local mappings (such that local derivations, $2$-local derivations, bilocal derivations, bilocal Lie derivations, weak-2-local derivations, local automorphisms,
$2$-local Lie $*$-automorphisms, 2-local $*$-Lie isomorphisms and so on) of associative algebras (see, for example, \cite{Nur}, \cite{aa17}, \cite{AEK21}, \cite{AEK211}, \cite{AK2016A}, \cite{AKO19}, \cite{Khrypchenko19}).
The study of local and $2$-local derivations of non-associative algebras was initiated in some papers of Ayupov and Kudaybergenov (for the case of Lie algebras, see~\cite{AK2016}, \cite{AK2016A}).
In particular, they proved that there are no non-trivial local and $2$-local derivations on complex semisimple finite-dimensional Lie algebras.
In~\cite{AK2016A} it is also given examples of $2$-local derivations on nilpotent Lie algebras which are not derivations.
After the cited works, the study of local and $2$-local derivations was continued for
many types of algebras, such that
Leibniz algebras~\cite{AKO19}, Jordan algebras~\cite{aa17},
$n$-ary algebras~\cite{FKK} and so on.
The first example of a simple (ternary) algebra with non-trivial local derivations is constructed by
  Ferreira, Kaygorodov and Kudaybergenov in~\cite{FKK}.
  After that, the first example of a simple (binary) algebra with non-trivial local derivations/automorphisms was constructed by
 Ayupov, Elduque and Kudaybergenov in~\cite{AEK21},\cite{AEK211}.
 The present paper is devoted to the study of local derivations of semisimple finite-dimensional Leibniz algebras.
 \begin{definition}
Let $\A$ be an algebra.
A linear map $\Delta : \A \rightarrow \A$ is called a local derivation, if for
any element $x \in \A$ there exists a derivation ${\mathfrak D}_x: \A \rightarrow \A$ such that $\Delta(x) = {\mathfrak D}_x(x)$.
\end{definition}

\section{Structure of semisimple Leibniz algebras and their derivations}

\subsection{Leibniz algebras}

Leibniz algebras present a "non antisymmetric" generalization of Lie algebras.
It appeared in some papers of Bloh [in 1960s]
and Loday [in 1990s].
Recently, they appeared in many geometric and physics applications (see, for example, \cite{bonez}, \cite{vert1}, \cite{leib2}, \cite{kotov20}, \cite{strow20} and references therein).
A systematic study of algebraic properties of Leibniz algebras is started from the Loday paper~\cite{L93}.
So, several classical theorems from Lie algebras theory have been extended
to the Leibniz algebras case;
many classification results regarding nilpotent, solvable, simple, and semisimple Leibniz algebras
are obtained
(see, for example, \cite{ak21}, \cite{AKO19}, \cite{AKOZ}, \cite{AOR19}, \cite{Barnes}, \cite{bm20}, \cite{ck20}, \cite{dzhuma}, \cite{FW21}, \cite{ls19}, \cite{mw21}, \cite{RMO}, \cite{t21}, and references therein).
Leibniz algebras is a particular case of terminal algebras and, on the other hand,
symmetric Leibniz algebras are Poisson admissible algebras.

An algebra $(\mathcal{L},[\cdot,\cdot])$ over a field $\mathbb{F}$ is called a
 (right)	Leibniz algebra if it satisfies the property
    \[
    [x,[y,z]] = [[x,y],z] - [[x,z],y],
    \]
	which is called Leibniz identity.
	For a Leibniz algebra $\mathcal{L}$, a subspace generated by its squares $\mathcal{I} = \text{span}\left\{[x,x]: x\in \mathcal{L}\right\}$
due to Leibniz identity becomes an ideal, and the quotient $\mathcal{G}_\mathcal{L} = \mathcal{L}/\mathcal{I}$ is a Lie algebra
called liezation of $\mathcal{L}.$ Moreover, $[ \mathcal{L}, \mathcal{I}] = 0.$
Following ideas of Dzhumadildaev~\cite{dzhuma},
	a Leibniz algebra $\mathcal{L}$ is called simple if its liezation is a simple Lie algebra and the ideal $\mathcal{I}$ is a
	simple ideal. Equivalently, $\mathcal{L}$ is simple iff $\mathcal{I}$ is the only non-trivial ideal of $\mathcal{L}.$
	A Leibniz algebra $\mathcal{L}$ is called {semisimple} if
	its liezation $\mathcal{G}_\mathcal{L}$ is a semisimple Lie algebra.
Simple and semisimple Leibniz algebras are under certain interest now \cite{AKO19}, \cite{AKOZ}, \cite{vert1}, \cite{dzhuma}, \cite{FW21}, \cite{mw21}, \cite{RMO}.

 Let $\mathcal{G}$ be a Lie algebra and $\mathcal{V}$ a (right) $\mathcal{G}$-module. Endow on vector space $\mathcal{L} = \mathcal{G} \oplus \mathcal{V}$ the bracket product as follows:
\[
	[(g_1, v_1), (g_2, v_2)]: = ([g_1, g_2], v_1 \cdot g_2),
\]
where $v \cdot g$ (sometimes denoted as $[v, g]$) is an action of an element $g$ of $\mathcal{G}$ on $v\in \mathcal{V}.$ Then
$\mathcal{L}$ is a Leibniz algebra, denoted as $\mathcal{G}\ltimes \mathcal{V}$. The following theorem proved by Barnes~\cite{Barnes} presents an analog of Levi-Malcev's theorem for Leibniz algebras.

\begin{theorem}
If $\mathcal{L}$ is a finite-dimensional Leibniz algebra over a field of characteristic zero, then $\mathcal{L} = \mathcal{S}\ltimes \mathcal{I},$ where $\mathcal{S}$ is a semisimple Lie subalgebra of $\mathcal{L}$.
\end{theorem}

It should be noted that $\mathcal{I}$ is a non-trivial module over the Lie algebra $\mathcal{S}$.
We say that a semisimple Leibniz algebra $\mathcal{L}$ is decomposable, if $\mathcal{L} = \mathcal{L}_1\oplus \mathcal{L}_2,$ where $\mathcal{L}_1$ and $\mathcal{L}_2$ are non-trivial semisimple Leibniz algebras.
Otherwise, we say that $\mathcal{L}$ is indecomposable.
Now we recall the structure of semisimple Leibniz algebras (see~\cite{AKOZ}).
Any complex semisimple finite-dimensional Leibniz algebra $\mathcal{L}$ represented as
\begin{align}\label{semisimplest}
\mathcal{L} = \bigoplus\limits_{i = 1}^n \left(\mathcal{S}_i\ltimes\mathcal{I}_i\right),
\end{align}
where each $\mathcal{S}_i\ltimes\mathcal{I}_i$ is an
indecomposable Leibniz algebra (see~\cite[Lemma 3.2]{AKOZ}).

\subsection{Derivations of semisimple Leibniz algebras}

Let $\mathcal{L}$ be a semisimple Leibniz algebra of the form~\eqref{semisimplest}.
It is clear that
\begin{align*}
\mathfrak{Der}\left(\mathcal{L}\right) = \bigoplus\limits_{i = 1}^n \mathfrak{Der}\left(\mathcal{S}_i\ltimes\mathcal{I}_i\right).
\end{align*}
Hence
\begin{align*}
\mathfrak{LDer}\left(\mathcal{L}\right) = \bigoplus\limits_{i = 1}^n \mathfrak{LDer}\left(\mathcal{S}_i\ltimes\mathcal{I}_i\right).
\end{align*}

So, it suffices to consider local derivations indecomposable semisimple Leibniz algebras.

Any indecomposable semisimple Leibniz algebra $\mathcal{L}$ represented as
\[
\mathcal{L} = \mathcal{S}\ltimes\left(\bigoplus\limits_{k = 1}^m \mathcal{V}_k\right),
\]
where each $\mathcal{S}\ltimes\mathcal{V}_k$ is also indecomposable semisimple Leibniz algebra.

Let $\mathcal{S}\ltimes\left(\bigoplus\limits_{k = 1}^m \mathcal{V}_k\right)$ be an indecomposable semisimple Leibniz algebra.
Denote by $\Gamma_\mathcal{S}$ the set of all $k = 1, \ldots, m$ such that $\mathcal{S}$ and $\mathcal{V}_k$ are isomorphic as $\mathcal{S}$-modules
and denote by $\Gamma_\mathcal{V}$ the set of all pairs $\{i,j\}$ such that $\mathcal{V}_i$ and $\mathcal{V}_j$ are isomorphic as $\mathcal{S}$-modules.

Any $\mathfrak{D}\in \mathfrak{Der}\left(\mathcal{S}\ltimes\left(\bigoplus\limits_{k = 1}^m \mathcal{V}_k\right)\right)$ is of the form
\begin{align}\label{derirred}
\mathfrak{D} = R_a+\sum\limits_{k\in \Gamma_\mathcal{S}}\varpi_k\theta_k+
\sum\limits_{\{i,j\}\in \Gamma_\mathcal{V}}\lambda_{i,j}\pi_{i,j},
\end{align}
where $\pi_{i,j}\in \textrm{Hom}_\mathcal{S}\left(\mathcal{V}_i, \mathcal{V}_j\right)$, $\theta_k\in \textrm{Hom}_\mathcal{S}\left(\mathcal{S}, \mathcal{V}_k\right)$ and $R_a$ is the standard right multiplication on $a \in \mathcal{S}$ (see~\cite[Theorem 4.5]{AKOZ}).

\section{Local derivations on semisimple Leibniz algebras}

The present part of the paper is dedicated to the proof of the following theorem.

\begin{theorem}\label{theoremlocal}
Let $\mathcal{L} = \mathcal{S}\ltimes \left(\bigoplus\limits_{k = 1}^m \mathcal{V}_k\right)$ be a complex semisimple finite-dimensional Leibniz algebra. Then any local derivation $\Delta$ on $\mathcal{L}$ is a derivation.
\end{theorem}

As we have mentioned above it suffices to consider local derivations of indecomposable semisimple Leibniz algebras.
From now on $\mathcal{L} = \mathcal{S}\ltimes \left(\bigoplus\limits_{k = 1}^m \mathcal{V}_k\right)$ is a
complex finite-dimensional indecomposable semisimple Leibniz algebra.

Let $\mathcal{H}$ be a Cartan subalgebra of $\mathcal{S}.$
Consider a root space decomposition of $\mathcal{S}:$
\[
\mathcal{S} = \mathcal{H}\oplus \bigoplus_{\alpha \in \Gamma}\mathcal{S}_\alpha,
\]
where $\Gamma$ is the set of all nonzero linear functionals $\alpha$ of $\mathcal{H}$ such that
\begin{center}
$\mathcal{S}_\alpha = \left\{x\in \mathcal{S}: [h, x] = \alpha(h)x,\, \forall\, h\in \mathcal{H}\right\}\neq \{0\}.$
\end{center}

Let
\[
\mathcal{V}_k = \bigoplus_{\beta \in \Phi_k}\mathcal{V}^\beta_k
\]
be a weight decomposition of $\mathcal{V}_k,$ where
$\Phi_k$ is the set of all weights.

For $q = 1, \ldots, m$ denote by $\textrm{pr}_q: \bigoplus\limits_{k = 1}^m \mathcal{V}_k\to \mathcal{V}_q$ a projection mapping defined as follows
\[
\textrm{pr}_q\left(\sum\limits_{k = 1}^m v_k\right) = v_q.
\]

Let us define a mapping $\Delta_{p,q}: \mathcal{S}\ltimes \mathcal{V}_p\to
\mathcal{V}_q$ as follows
\[
\Delta_{p,q}(x+v) = \textrm{pr}_q\left(\Delta(x+v)\right),\,\, x+v\in \mathcal{S}\ltimes \mathcal{V}_p.
\]
By~\eqref{derirred} for any $x+v\in \mathcal{S}\ltimes \mathcal{V}_p^\beta$ there exist $a_{x+v}\in \mathcal{S}$ and complex numbers $\omega_k^{(x+v)}, \lambda_{i,j}^{(x+v)}$ such that
\begin{align*}
\Delta(x+v) = [x+v, a_{x+v}]+\sum\limits_{k\in \Gamma_\mathcal{S}}\omega_k^{(x+v)}\theta_k(x)+\sum\limits_{\{i,j\}\in \Gamma_\mathcal{V}}\lambda_{i,j}^{(x+v)}\pi_{i,j}(v).
\end{align*}
Let $a_{x+v} = h_{x+v}+\sum \limits_{\alpha \in \Gamma} c_\alpha^{(x+v)} e_\alpha\in \mathcal{H}\oplus \oplus_{\alpha\in \Gamma} \mathcal{S}_\alpha$, and denote $\Gamma_q = \left\{\alpha: [v, e_\alpha]\in \mathcal{V}_q\right\}$. Then
\begin{align*}
\Delta_{p,q}(x+v) & = & \textrm{pr}_q\left([x+v,a_{x+v}]+\sum\limits_{k\in \Gamma_\mathcal{S}}\omega_k^{(x+v)}\theta_k(x)
+\sum\limits_{\{i,j\}\in \Gamma_\mathcal{V}}\lambda_{i,j}^{(x+v)}\pi_{i,j}(v)\right),
\end{align*}
that is
\begin{align}\label{derres}
\Delta_{p,q}(x+v) & = & \delta_{p,q} [v,h_{x+v}]+
\left[v, \sum\limits_{\alpha\in \Gamma_q}c_\alpha^{(x+v)} e_\alpha\right]+\omega_q^{(x+v)}\theta_q(x)+
\lambda_{p,q}^{(x+v)}\pi_{p,q}(v),
\end{align}

 where $\delta_{p,q}$ is the Kronecker delta.

If necessary, after renumbering we can assume that there is a number $k\in \{1, \ldots, m\}$ such that
$\mathcal{S}$ and $\mathcal{V}_i$ are isomorphic
\(\mathcal{S}\)-modules for all $i = 1, \ldots, k,$ and
$\mathcal{S}$ and $\mathcal{V}_i$ are not isomorphic
\(\mathcal{S}\)-modules for all $i = k+1, \ldots, m.$
According~\eqref{derirred}, we can consider two possible cases separately.

Let us first consider an indecomposable semisimple Leibniz algebra
$\mathcal{L} = \mathcal{S}\ltimes \left(\bigoplus\limits_{i = 1}^m \mathcal{V}_i\right)$ such that
\(\mathcal{S}\) and \(\mathcal{V}_i\, (i = 1, \ldots, m)\) are isomorphic
\(\mathcal{S}\)-modules.

\begin{lemma}\label{GeqI}
Let \(\Delta\) be a local derivation on \(\mathcal{L}\) such that
\(\Delta (\mathcal{L}) \subseteq \mathcal{V}\). Then
\(\Delta\) is a derivation.
\end{lemma}

\begin{proof}

Fix the indices $p, q.$
Let us show that there are complex numbers $\omega_q$ and $\lambda_{p,q}$ such that
\[
\Delta_{p,q} = \omega_q\theta_q+\lambda_{p,q}\pi_{p,q}.
\]
Fix a basis \(\{x_1, \ldots, x_m\}\) in \(\mathcal{S}\).
The system of vectors \(\left\{\theta_s(x_i)\right\}_{1 \leq i\leq m}\) is a basis in \(\mathcal{V}_s\) \ (for $s = p,q$).
Here,
\(\theta_s\) is an \(\mathcal{S}\)-module isomorphism from
\(\mathcal{S}\) onto \(\mathcal{V}_s\), in particular,
\[
    \theta_s([x,y]) = [\theta_s(x), y].
\]

Using~\eqref{derres} for \(x = x_i \) and take a complex number \(\omega_{q}^{(i)}\)
such that
\[
\Delta_{p,q}(x_i) = \omega_{q}^{(i)}\theta_q(x_i).
\]

Now for the element \(x = x_i+x_j\), where \(i\neq j\), take a complex number \(\omega_{q}^{(i,j)} \)
such that
\[
\Delta_{p,q}(x_i+x_j) = \omega_{q}^{(i,j)}\theta_q(x_i+x_j) = \omega_{q}^{(i,j)}\theta_q(x_i)+\omega_{q}^{(i,j)}\theta_q(x_j).
\]
On the other hand,
\[
\Delta_{p,q}(x_i+x_j) = \omega_{q}^{(i)}\theta_q(x_i)+\omega_{q}^{(j)}\theta_q(x_j).
\]
Comparing the last two equalities we obtain \(\omega_{q}^{(i)} = \omega_{q}^{(j)}\) for all \(i, j\). This means that there
exists a complex number \(\varpi_q \) such that
\begin{align}\label{xxx}
\Delta_{p,q}(x_i) = \varpi_q \theta_q(x_i).
\end{align}
Now by~\eqref{derres} for \(x = x_i+\theta_p(x_i)\in \mathcal{S}\ltimes \mathcal{V}_p\) take complex numbers \(\omega_i\) and \(\lambda_i\) such that
\[
\Delta(x_i+\theta_p(x_i)) = \omega_i\theta_q(x_i)+\lambda_i
\pi_{p,q}(\theta_p(x_i)) = \left(\omega_i+\lambda_i
\right)\theta_q(x_i).
\]
Taking into account (\ref{xxx}) we obtain that
\begin{longtable}
{lclcl}
$\Delta_{p,q}(\theta_p(x_i))$ &$ = $&$\Delta_{p,q}(x_i+\theta_p(x_i))-\Delta_{p,q}(x_i)$& \\
 & $ = $ & $(\omega_i+\lambda_i)\theta_q(x_i)
-\varpi_q\theta_q(x_i)$&$ = $&$(\omega_i-\varpi_q+\lambda_i)\theta_q(x_i).$
\end{longtable}
This means that for every \(i\in \{1, \ldots, m\}\) there exists a complex
number \(\Lambda_i\) such that
 \begin{align}\label{ii}
\Delta(\theta_p(x_i)) = \Lambda_i \theta_q(x_i).
 \end{align}
Take an element \(x = x_i+x_j+\theta_p(x_i+x_j)\in
\mathcal{S}\ltimes\mathcal{V}_p\), where \(i\neq j\). By~\eqref{derres}, we get that
\[
\Delta_{p,q}(x_i+x_j+\theta_p(x_i+x_j)) = \omega_{i,j}\theta_q(x_i+x_j)+
\lambda_{i,j}\theta_q(x_i+x_j).
\]
Taking into account (\ref{xxx}) we obtain that
\begin{align*}
\Delta_{p,q}\left(\theta_q(x_i+x_j)\right)
& = \Delta_{p,q}(x_i+x_j+\theta_q(x_i+x_j))-\Delta_{p,q}(x_i+x_j) \\
& = (\omega_{i,j}-\varpi_q+\lambda_{i,j})\theta_q(x_i+x_j).
\end{align*}
On the other hand, by~\eqref{ii},
\[
\Delta_{p,q}(\theta_q(x_i+x_j)) = \Delta_{p,q}(\theta_q(x_i))+\Delta_{p,q}(\theta_q(x_j))
= \lambda_i \theta_q(x_i)+\lambda_j\theta_q(x_j).
\]
Comparing the last two equalities we obtain that
\(\lambda_i = \lambda_j\) for all \(i\) and \(j\). This means that
there exist a complex number \(\lambda_{p,q}\) such that
\begin{align}\label{yyy}
\Delta_{p,q}(\theta_q(x_i)) = \lambda_{p,q}\theta_q(x_i).
\end{align}
Combining (\ref{xxx}) and (\ref{yyy})
we obtain that \(\Delta_{p,q} = \varpi_q\theta_q+\lambda_{p,q}\pi_{p,q}\). This means that
\(\Delta\) is a derivation. The proof is completed.
\end{proof}

In the next lemma we consider \(\mathcal{L} = \mathcal{S}\ltimes \left(\bigoplus\limits_{k = 1}^m \mathcal{V}_k\right)\), an indecomposable semisimple Leibniz algebra, such that \(\mathcal{S}\) and \(\mathcal{V}_k\) are not isomorphic \(\mathcal{S}\)-modules for all $k = 1, \ldots, m.$

\begin{lemma}\label{GeqII}
Let \(\Delta\) be a local derivation on \(\mathcal{L}\) such that
\(\Delta (\mathcal{L}) \subseteq \mathcal{V}\). Then
\(\Delta\) is a derivation.
\end{lemma}
\begin{proof}
Let \(\left\{v_1^{(1)}, \ldots, v_n^{(1)}\right\}\) be a basis of \(\mathcal{V}_1\). Since $\mathcal{V}_1$ and $\mathcal{V}_k$ are isomorphic, it follows that \(\left\{v_i^{(q)} = \pi_{1,q}(v_i^{(1)}): i = 1,\ldots n\right\}\) is a basis of \(\mathcal{V}_q\).

Without lost of generality we can
assume that for any \(v_i^{(1)}\) there exists a weight~\(\beta_i\) such
that \(v_i^{(1)}\in \mathcal{V}_{\beta_i}\).
Let \(h_0\) be a strongly regular element in \(\mathcal{H}\), that is, $\alpha(h_0)\neq \beta(h_0)$ for any $\alpha, \beta\in \Gamma, \alpha\neq \beta.$ For
\(x = h_0+v_i^{(1)}\in \mathcal{S}\ltimes\mathcal{V}_1\) take an element \(a_x\in
\mathcal{S}\) and complex numbers \(\lambda^{(x)}_k \) such that
\[
\Delta\left(h_0+v_i^{(1)}\right) = [h_0, a_x]+\left[v_i^{(1)}, a_x\right]+\sum\limits_{k\in \Gamma_{1,k}}\lambda^{(x)}_k v_i^{(k)}.
\]
Taking into account that $h_0$ is strongly regular, from \([h_0, a_x] = 0\), we have that \(a_x\in \mathcal{H}\).
Further
\begin{align*}
\Delta\left(v_i^{(1)}\right)
&= \Delta\left(h_0+v_i^{(1)}\right) \\
&= \left[v_i^{(1)}, a_x\right]+\sum\limits_{k\in \Gamma_{1,k}} \lambda^{(x)}_k v_i^{(k)} = (\beta_i(a_x)+\lambda^{(x)}_1)v_i^{(1)}+\sum\limits_{1<k\in \Gamma_{1,k}}\lambda^{(x)}_k v_i^{(k)}.
\end{align*}
Now we change the element \(x = h_0+v_i^{(1)}\) to the element \(\overline{x} = h_0+v_i^{(1)}+v_j^{(1)}\) $(i\neq j),$ then similar as above we get that
\[
\Delta\left(v_i^{(1)}+v_j^{(1)}\right) = (*)v_i^{(1)}+(*)v_j^{(1)}+\sum\limits_{1<k\in \Gamma_{1,k}}\lambda^{(\overline{x})}_k \left(v_i^{(k)}+v_j^{(k)}\right).
\]
Comparing the last two equalities we can see that
there are $\lambda_2, \ldots, \lambda_n\in \mathbb{C}$ such that
\[
\Delta\left(v_i^{(1)}\right) = (*)v_i^{(1)}+\sum\limits_{1<k\in \Gamma_{1,k}}\lambda_k v_i^{(k)}.
\]
Replacing $\Delta$ with $\Delta-\sum\limits_{1<k\in \Gamma_{1,k}}
\lambda_k \pi_{1,k}$ we obtain a new local derivation which maps $\mathcal{V}_1$ into itself.
Due to~\eqref{derirred}, there exist complex numbers \(\lambda_i, i = 1, \ldots, n\) such that
\begin{align}\label{res}
\Delta\left(v_i^{(1)}\right) = \lambda_i v_i^{(1)}.
\end{align}
We shall show that \(\lambda_1 = \ldots = \lambda_n\). For a fixed \(v_i^{(1)}\) \((i\neq 1)\) we have that
\begin{align}\label{lambda}
\Delta\left(v_1^{(1)}+v_i^{(1)}\right) = \lambda_1v_1^{(1)}+
\lambda_iv_i^{(1)}.
\end{align}

Without loss of generality we can assume that \(\beta_1\) is a
fixed highest weight of \(\mathcal{V}_1\). It is known~\cite[Page
108]{Hum} that the weight \(\beta_1-\beta_i\) can be represented as
\[
\beta_1-\beta_i = n_1\alpha_1+\ldots+n_l\alpha_l,
\]
where \(\alpha_1, \ldots, \alpha_l\) are simple roots of
\(\mathcal{S}\), \(n_1, \ldots, n_l\) are non negative integers.

Below we shall consider two separated cases.
\begin{enumerate}
 \item[\textbf{Case 1.}]
\(\alpha_0 = n_1\alpha_1+\ldots+n_l\alpha_l\) is not a root.
Take the following element
\[
y = n_1e_{\alpha_1}+\ldots+n_le_{\alpha_l}+v_1^{(1)}+v_i^{(1)}.
\]
By the definition of local derivation we can find an element \(a_y = h+\sum\limits_{\alpha\in \Gamma}c_\alpha
e_\alpha\in \mathcal{S}\) and a number \(\lambda^{(y)}\) such that
\[
\Delta(y) = [y, a_y]+\lambda^{(y)}\left(v_1^{(1)} +v_i^{(1)}\right).
\]
Taking into account~\eqref{res} and \(\Delta(y)\in \mathcal{V}\), we obtain that
\[
\left[\sum\limits_{s = 1}^l n_s e_{\alpha_s},
h+\sum\limits_{\alpha\in \Gamma}c_\alpha e_\alpha\right] = 0.
\]
Thus
\[
\sum\limits_{s = 1}^l n_s
\alpha_s(h)e_{\alpha_s}+\sum\limits_{t = 1}^l \sum\limits_{\alpha\in
\Gamma} (\ast) e_{\alpha+\alpha_t} = 0,
\]
where the symbols \((\ast)\) denote appropriate coefficients. The
second summand does not contain any element of the form
\(e_{\alpha_s}\). Indeed, if we assume that
\(\alpha_s = \alpha+\alpha_t\), we have that
\(\alpha = \alpha_s-\alpha_t\). But $\alpha_s-\alpha_t$ is not a
root, because \(\alpha_s, \alpha_t\) are simple roots. Hence all
coefficients of the first summand are zero, i.e.,
\[
n_1 \alpha_1(h) = \ldots = n_l \alpha_l(h) = 0.
\]
Further
\[
\Delta\left(v_1^{(1)} +v_i^{(1)}\right) = \Delta(y) = \left[v_1^{(1)} +v_i^{(1)}, a_x\right]+\lambda^{(y)}\left(v_1^{(1)} +v_i^{(1)}\right).
\]
Let us calculate the product \(\left[v_1^{(1)} +v_i^{(1)}, a_x\right]\).
We have
\begin{align*}
\left[v_1^{(1)} +v_i^{(1)}, a_x\right] 
& = \left[v_1^{(1)} +v_i^{(1)}, h+\sum\limits_{\alpha\in
\Phi}c_\alpha e_\alpha\right]\\
& = \beta_1(h)v^{(1)}_{\beta_1}+\beta_2(h)v^{(1)}_{\beta_2}+\sum\limits_{t = 1}^2 \sum\limits_{\alpha\in
\Phi}(\ast) v^{(1)}_{\beta_t+\alpha}.
\end{align*}
The last summand does not contain \(v^{(1)}_{\beta_1}\) and
\(v^{(1)}_{\beta_i}\), because \(\beta_1-\beta_i\) is not a root by the assumption. This means that
\begin{align}\label{iii}
\Delta\left(v^{(1)}_{\beta_1}+v^{(1)}_{\beta_i}\right) = \left(\beta_1(h)+\lambda^{(y)}\right)v^{(1)}_{\beta_1}+\left(\beta_i(h)+\lambda^{(y)}\right)v^{(1)}_{\beta_i}.
\end{align}
The difference of the coefficients of the right side is
\[
\beta_1(h)-\beta_i(h) = \sum\limits_{s = 1}^l n_s\alpha_s(h) = 0,
\]
because of \(n_1\alpha_1(h) = \ldots = n_l\alpha_l(h) = 0\). Finally,
comparing coefficients in (\ref{lambda}) and (\ref{iii}) we get
\[
\lambda_1 = \beta_1(h)+\lambda^{(y)} = \beta_i(h)+\lambda^{(y)} = \lambda_i.
\]

\item[\textbf{Case 2.}] \(\alpha_0 = n_1\alpha_1+\ldots+n_l\alpha_l\) is a root.
Note that \(\dim
\mathcal{V}_{\beta_1} = 1\), because \(\beta_1\) is a highest weight. Since \(\beta_1-\beta_i\) is
a root,~\cite[Lemma 3.2.9]{Good} implies that \(\dim
\mathcal{V}_{1, \beta_i} = \dim \mathcal{V}_{1, \beta_1}\), and hence
there exist numbers \(t_{-\alpha_0}\neq0\) and \(t_{\alpha_0}\)
such that
\[
\left[v^{(1)}_{\beta_1},
e_{-\alpha_0}\right] = t_{-\alpha_0}v^{(1)}_{\beta_i},\, \left[v^{(1)}_{\beta_i},
e_{\alpha_0}\right] = t_{\alpha_0}v^{(1)}_{\beta_1}.
\]
Take the following element
\[
z = t_{-\alpha_0}e_{\alpha_0}+t_{\alpha_0}e_{-\alpha_0}+v^{(1)}_{\beta_1}+v^{(1)}_{\beta_i},
\]
and choose an element \(a_z = h+\sum\limits_{\alpha\in \Phi}c_\alpha
e_\alpha\in \mathcal{S}\) and a number \(\lambda_z\) such that
\[
\Delta(z) = [z, a_z]+\lambda_z\left(v^{(1)}_{\beta_1}+v^{(1)}_{\beta_i}\right).
\]
Since \(\Delta(z)\in \mathcal{V}\), we obtain that
\[
\left[t_{-\alpha_0}e_{\alpha_0}+t_{\alpha_0}e_{-\alpha_0},
h+\sum\limits_{\alpha\in \Phi}c_\alpha e_\alpha\right] = 0.
\]
Now rewrite the last equality as
\[
\alpha_0(h)t_{-\alpha_0}e_{\alpha_0}-\alpha_0(h)t_{\alpha_0}e_{-\alpha_0}+(t_{-\alpha_0}c_{-\alpha_0}-
t_{\alpha_0}c_{\alpha_0})h_{\alpha_0}+\sum\limits_{\alpha\neq\pm\alpha_0}
(\ast) e_{\alpha\pm\alpha_0} = 0,
\]
where \(h_{\alpha_0} = [e_{\alpha_0}, e_{-\alpha_0}]\in
\mathcal{H}\). The last summand in the sum does not contain
elements \(e_{\alpha_0}\) and \(e_{-\alpha_0}\). Indeed, if we
assume that \(\alpha_0 = \alpha-\alpha_0\), we have that
\(\alpha = 2\alpha_0\). But \(2\alpha_0\) is not a root. Hence the
first three coefficients of this sum are zero, i.e.,
\begin{align}\label{zzz}
\alpha_0(h) = 0,\,
t_{\alpha_0}c_{\alpha_0} = t_{-\alpha_0}c_{-\alpha_0}.
\end{align}
Further
\[
\Delta\left(v^{(1)}_{\beta_1}+v^{(1)}_{\beta_i}\right) = \Delta(z) = \left[v^{(1)}_{\beta_1}+v^{(1)}_{\beta_i},
a_z\right]+\lambda_z\left(v^{(1)}_{\beta_1}+v^{(1)}_{\beta_i}\right).
\]
Let us consider the element \(\left[v^{(1)}_{\beta_1}+v^{(1)}_{\beta_i}, a_z\right]\).
We have
\begin{align*}
\left[v^{(1)}_{\beta_1}+v^{(1)}_{\beta_i}, a_z\right]
&= \left[v^{(1)}_{\beta_1}+v^{(1)}_{\beta_i}, h+\sum\limits_{\alpha\in
\Phi}c_\alpha e_\alpha\right] \\
&= \left[v^{(1)}_{\beta_1}, h\right]
    +c_{\alpha_0}\left[v^{(1)}_{\beta_i}, e_{\alpha_0}\right]
    +\left[v^{(1)}_{\beta_i}, h\right] \\
&{\quad}+c_{-\alpha_0}\left[v^{(1)}_{\beta_1}, e_{-\alpha_0}\right]
    + c_{\alpha_0}\left[v^{(1)}_{\beta_1}, e_{\alpha_0}\right]
    + c_{-\alpha_0}\left[v^{(1)}_{\beta_i}, e_{-\alpha_0}\right] \\
&{\quad}+\sum\limits_{\alpha\neq\pm\alpha_0}c_\alpha\left[v^{(1)}_{\beta_1}, e_\alpha\right]
    + \sum\limits_{\alpha\neq\pm\alpha_0}c_\alpha\left[v^{(1)}_{\beta_i}, e_\alpha\right] \\
&= (\beta_1(h)+t_{\alpha_0}c_{\alpha_0})v^{(1)}_{\beta_1}
    +(\beta_i(h)+t_{-\alpha_0}c_{-\alpha_0})v^{(1)}_{\beta_i} \\
&{\quad}+(\ast) v^{(1)}_{2\beta_1-\beta_i}
    +(\ast) v^{(1)}_{2\beta_i-\beta_1} \\
&{\quad}+ \sum\limits_{\alpha\neq\pm\alpha_0}(\ast) v^{(1)}_{\beta_1+\alpha}
    + \sum\limits_{\alpha\neq\pm\alpha_0}(\ast) v^{(1)}_{\beta_i+\alpha}.
\end{align*}
The last three summands do not contain \(v^{(1)}_{\beta_1}\) and \(v^{(1)}_{\beta_i}\), because \(\beta_1-\beta_i = \alpha_0\) and
\(\alpha\neq\pm\alpha_0\).
This means that
\begin{align}\label{ttt}
\Delta\left(v^{(1)}_{\beta_1} + v^{(1)}_{\beta_i}\right)
= {}&{}(\beta_1(h) + t_{\alpha_0}c_{\alpha_0} + \lambda_z)v^{(1)}_{\beta_1} \\
&{}+ (\beta_i(h)+t_{-\alpha_0}c_{-\alpha_0} + \lambda_z)v^{(1)}_{\beta_i}. \nonumber
\end{align}
Taking into account (\ref{zzz}) we find the difference of coefficients on the right side:
\[
(\beta_1(h)+t_{\alpha_0}c_{\alpha_0})-(\beta_i(h)+t_{-\alpha_0}c_{-\alpha_0}) =
\alpha_0(h)+t_{\alpha_0}c_{\alpha_0}-t_{-\alpha_0}c_{-\alpha_0} = 0.
\]
Combining (\ref{lambda}) and (\ref{ttt}) we obtain that
\[
\lambda_1 = \beta_1(h)+t_{\alpha_0}c_{\alpha_0}+\lambda_z = \beta_i(h)+t_{-\alpha_0}c_{-\alpha_0}+\lambda_z = \lambda_i.
\]
So, we have proved that
\(\Delta\left(v_i^{(1)}\right) = \lambda_1 v^{(1)}_i\)
for all $i = 1,\ldots, n.$ By a similar way we obtain that
\(\Delta\left(v_i^{(k)}\right) = \lambda_k v^{(k)}_i\)
for all $i = 1,\ldots, n_k.$ Thus \(\Delta = \sum\limits_{k = 1}^m \lambda_k \pi_{k,k}\), and therefore
\(\Delta\) is a derivation.
\end{enumerate}
The proof is completed.
\end{proof}

\begin{proof}[\textbf{Proof of Theorem~\ref{theoremlocal}}.]
Let \(\Delta\) be an arbitrary local derivation on
\(\mathcal{L}\). For an arbitrary element \(x\in \mathcal{S}\)
take a derivation $\mathfrak{D} = R_{a_x}+\sum\limits_{k\in \Gamma_\mathcal{S}}\varpi_k^{(x)}\theta_k^{(x)}+
\sum\limits_{\{i,j\}\in \Gamma_\mathcal{V}}\lambda_{i,j}^{(x)}\pi_{i,j}^{(x)}$ of the form~\eqref{derirred} such that
\[
\Delta(x) = [x, a_x]+\sum\limits_{k\in \Gamma_\mathcal{S}}\omega_{x,k}^{(x)}\theta_k^{(x)}(x).
\]
Then the mapping
\[
x\in \mathcal{S}\rightarrow [x, a_x]\in \mathcal{S}
\]
is a well-defined local derivation on \(\mathcal{S}\), and
by~\cite[Theorem 3.1]{AK2016} it is a derivation
generated by an element \(a\in \mathcal{S}\). Then the local
derivation \(\Delta-R_{a_x}\) maps \(\mathcal{L}\) into
\(\mathcal{V}\). By Lemmas~\ref{GeqI} and~\ref{GeqII} we get that \(\Delta-R_{a_x}\) is
a derivation and therefore \(\Delta\) is also a derivation.
The proof is completed.
\end{proof}

\EditInfo{%
    January 03, 2022}{%
    January 28, 2022}{%
    Friedrich Wagemann
    }

\end{paper}